\documentclass[12pt]{amsart}
\usepackage{latexsym}

\usepackage{amsmath,amssymb,amsfonts,array}

\def\CC{\mathbb{C}}

\def\RR{\mathbb{R}}

\def\delbar{\overline{\partial}}

\address{Stanford University, Department of
  Mathematics, Stanford CA 94305--2125, USA, Supported by the DFG}
\email{mohnke@math.stanford.edu}
\title{Holomorphic Disks and the Chord Conjecture}
\author{Klaus Mohnke}
\date{}

\newtheorem{theorem}{Theorem}

\newtheorem{proposition}[theorem]{Proposition}

\theoremstyle{definition}

\theoremstyle{remark}

\theoremstyle{plain}
\newtheorem*{Fundamental-lem}{Fundamental Lemma}

\begin{document}
\bibliographystyle{plain}

\begin{abstract}
We prove ( a weak version of)
Arnold's Chord Conjecture in \cite{Arnold:firststeps} using
Gromov's ``classical'' idea in \cite{Gromov} to produce holomorphic
disks with boundary on a Lagrangian submanifold.
\end{abstract}

\maketitle

\subsection*{Arnold's Chord Conjecture}
In this paper we prove the following theorem which was
conjectured by Arnold \cite{Arnold:firststeps}:
\begin{theorem}\label{arnold}
For every closed, compact Legendrian submanifold in $S^{2n-1}$
with the standard contact structure and any contact form
for this structure, there is a Reeb chord, i.e.~an integral curve of the Reeb
vector field which begins and ends on the Legendrian submanifold.
\end{theorem}
Theorem~\ref{arnold} will follow as a corollary from the main result of this
paper, Theorem~\ref{main}. In fact it can be applied to a
more general situation:
\begin{theorem}\label{chord}
Let $(M,\xi)$ be a contact structure which arises as smooth
boundary of a compact subcritical Stein manifold (see 
\cite{Biran/Cieliebak:subcritical} for a definition).
Then for any Legendrian
submanifold  and any contact one form corresponding to $\xi$
there is a Reeb chord.
\end{theorem}
Our results include the existence of chords for Legendrians
in the standard contact structure on $\RR P^{2n-1}$ proved
by Ginzburg and Givental \cite{Givental:MaslovLMS, Givental:Maslov},
although it does not provide their statement of linear growth.
They cover results by Abbas \cite{Abbas:chord}
and Cieliebak \cite{Cieliebak:chord} who treat subcases of
the problem on the sphere and on boundaries of subcritical
Stein manifolds.

\subsection*{Lagrangian out of Legendrian embeddings}
Consider a closed Legendrian submanifold $l\subset M^{2n-1}$ in a
contact manifold. Given a contact one form $\alpha$
we will construct Lagrangian embeddings of the torus
$l\times S^1$ into the symplectization
$(M\times\RR, d(e^s\alpha))$ of $(M,\xi)$ and study them.
Denote by $\phi$ the flow of the Reeb vector field $R=R_{\alpha}$.

Assume $l$ has no Reeb chords of length at most $T>0$.
Then there is an embedding
$\Phi:l\times(\RR\times[0,T])\hookrightarrow M\times\RR$
such that $\Phi^{\ast}(d(e^{s}\alpha))=e^{s}ds\wedge dt$,
$(s,t)\in\RR\times[0,T]$ being the coordinates of 
the infinite strip. The map is, of course, constructed
using the Reeb flow: $\Phi(x,s,t):=(\phi^{t}(x),s)$.

The upshot of all this is, that any smooth embedded loop
in the rectangle $[S,0]\times[0,T]$ defines a smooth Lagrangian
embedding of $l\times S^{1}$ into $M\times [S,0]$.
Here $S<0$ is any negative parameter. The integral of the primitive
$e^{s}\alpha$ over closed curves in such a Lagrangian
depends only on the homology class of that curve.
It vanishes on classes arising as closed loops in $l$,
since $\alpha|_{l}\equiv 0$. On $\{x\}\times l$ it is given
by the area enclosed by the loop in the rectangle
with respect to the volume form $e^{s}dsdt$.
Approximating the boundary of $[S,0]\times[0,T]$
by smooth embedded loops we obtain the following
\begin{proposition}
Assume the Legendrian $l\subset M$ in a contact manifold
has no Reeb chord of length smaller or equal to $T$ with respect to
a contact one form $\alpha$. Then
there exists a smooth Lagrangian embedding
of $l\times S^1$ into the symplectic manifold
$(M\times [S,0], d(e^s\alpha))$ such that the
integral of the primitive $e^s\alpha$
over a smooth closed curve in the image is
an integer multiple of some constant $C>(1-e^S)T$.
\end{proposition}

\subsection*{Proof of Arnold's Chord Conjecture}
The results described above will follow from  the
following
\begin{theorem}\label{main}
Assume $W$ is a symplectic manifold of bounded
geometry, i.e.~there is an almost complex structure
$J$ together with a hermitian metric $g$ whose
sectional curvature is globally bounded from
above and injectivity radius bounded from below, such
that $\omega$ tames $J$ uniformly w.r.t.~$g$:
$
\omega(X,JX)\ge \mbox{const.}\|X\|^{2}_{g}
$
for any tangent vector $X\in TW$. Suppose $W$ is
monotone, i.e.~any sphere in $W$ has vanishing symplectic area.
Let $(M\times[S,0], d(e^s\alpha))\hookrightarrow (W,\omega)$
be a symplectic embedding of a finite cylinder in
the symplectization of a contact manifold into $W$, such
that $\pi_{1}(M)$ injects into $\pi_{1}(W)$. 
Suppose the image of the embedding can be
displaced by a (time-dependent) Hamiltonian flow with compact support
and oscillation
$$
\|H\|:= \int_{0}^{1}\Big(\max_{x\in W}H(t,x) - \min_{x\in 
W}H(t,x)\Big)dt.
$$
Then {\em any} closed Legendrian
$l\subset M$ admits a Reeb chord of length not bigger than
$\|H\|/(1-e^S)$.
\end{theorem}

\begin{proof}Assume that there is a Legendrian embedding into $M$ which
admits no Reeb chord of length not bigger than $T>0$.
Consider a Lagrangian embedding
as constructed in the previous section. It will be displaced
from itself by $H$.
Fix any almost complex structure $J$ as in the theorem.
Due to monotonicity of $W$ there are no $J$--holomorphic
spheres in $W$. Then an argument
given by  Chekanov in \cite{Chekanov:displacement} produces a
non--constant $J$--holomorphic disk in $W$
with boundary on $L$ and (symplectic) area  smaller than or equal to
$\|H\|$. But due to the
injectivity of $\pi_{1}(M)\longrightarrow\pi_{1}(W)$ there
is another disk completely lying in $M$. Since $W$ is monotone
the symplectic areas of both disks agree. The latter is an
integer multiple of $C>T(1-e^{S})$. Non-constant holomorphic
disks have positive symplectic area, hence $T\le \|H\|/(1-e^{S})$. 
\end{proof}

{\em Remarks.} (1) Arnold's conjecture follows from this result since
every contact one form corresponding to the standard
contact structure can be realized as the restriction of
the primitive $\omega(z,.)$ of
the standard symplectic structure in $\CC^{n}$ to the smooth boundary
of a star-shaped domain.\\
(2) Biran and Cieliebak  observed in \cite{Biran/Cieliebak:subcritical}
that {\em any} compact set in a (complete)
subcritical Stein manifold can be displaced from itself by a 
Hamiltonian isotopy with compact support. In particular, one can apply this to the
strongly pseudo--convex domain in the Stein manifold within
its completion to obtain Theorem~\ref{chord}.\\
(3) One could prove the result using Gromov's Fredholm alternative
for maps of the disk into $W$ with boundary on a Lagrangian
satisfying a $\delbar$--equation with a non-zero right-hand side
(see \cite{Gromov} and \cite{Audin/Lalonde/Polterovich}). It
should be an amusing and (non-trivial) exercise to produce exactly
the same estimate. For our convenience we chose to use Chekanov's
result but feel obliged to point out that the result is
a consequence of Gromov's original ideas.
\newline

{\em Problem.} We have not answered Arnold's
original question yet, whether there always exists a Reeb orbit
intersecting the Legendrian submanifold at least twice.
We proved that there is a chord whose length is
estimated by $S$ and the displacement energy. But this chord could
be a closed orbit which intersects the Legendrian once.
Notice that this question is bound to the sphere. E.g.~pick the
contact form $x_{1}d\theta + \frac{1}{2}(x_{2}dx_{3}-x_{3}dx_{2})$
on $S^{1}\times S^{2}$ with coordinate $\theta\in S^{1}$ and
$S^{2}=\{x=(x_{1},x_{2},x_{3})\in \RR^{3}\mid \|x\|=1\}$. Then
$S^{1}\times\{0\}$ is a Legendrian loop which intersects each of the
(closed) Reeb orbits $\{\theta\}\times \{x_{1}=0\}$ only once.
\newline

{\em Acknowledgements.} I would like to thank Leonid Polterovich,
Kai Cieliebak, Yasha Eliashberg, Tobias Ekholm and John Etnyre
for their interest and comments, and the
Department of Mathematics of Stanford for the warm hospitality during
my visit.

\end{document}